\documentclass[11pt]{article}
\usepackage{mathrsfs}
\usepackage{latexsym,lineno}
\usepackage{epsfig}
\usepackage{color}
\usepackage{amsmath}\usepackage{fleqn}\usepackage{verbatim}\usepackage{epsf}
\usepackage{amsthm}\usepackage{graphicx, float}\usepackage{graphicx}
\usepackage{amsfonts}\usepackage{amssymb}\usepackage{graphpap}
\usepackage{epic}\usepackage{curves}
\newtheorem{theorem}{Theorem}[section]

\newtheorem{proposition}{Proposition}[section]
\newtheorem{observation}{Observation}[section]
\newtheorem{definition}{Definition}[section]

\newtheorem{lemma}{Lemma}[section]

\title{\bf Matching criticality in intersecting hypergraphs\thanks {Research was partially supported by  the National Nature Science Foundation of China (Nos. 11571222, 11471210)}}
\author {Liying Kang$^{1}$, \, Zhenyu Ni$^{1}$, \,  Erfang Shan$^{1,2}$\thanks{\em Corresponding authors. Email address: efshan@shu.edu.cn (E. Shan)}\\
{\small $^{1}$Department of Mathematics, Shanghai University,
Shanghai 200444, P.R. China}\\
{\small$^{2}$School of Management, Shanghai University,
Shanghai 200444, P.R. China}}

\date{}
\begin{document}

\maketitle

\begin{abstract}
The transversal number $\tau(H)$ of a hypergraph $H$ is the minimum cardinality of a set
of vertices that intersects all edges of $H$.
The matching number $\alpha'(H)$ of $H$ is the size of a largest matching in $H$, where a matching is a set of pairwise
disjoint edges in $H$.
  A hypergraph is intersecting if each pair of edges has a
nonempty intersection. Equivalently, $H$ is an intersecting hypergraph if and only if $\alpha'(H)=1$.
We observe that $\tau(H)\le r$ for an intersecting hypergraph $H$ of rank $r$.
 For an intersecting hypergraph $H$ of rank $r$ without isolated vertex, we call $H$   1-special if $\tau(H)=r$;  $H$ is maximal 1-special if $H$ is 1-special and adding any missing $r$-edge to $H$ increases the matching number. Furthermore, $H$ is called 1-edge-critical if for any $e\in E(H)$ and any $v\in e$, $v$-shrinking $e$ increases the matching number; $H$ is called 1-vertex-critical if for every vertex $v$ in $H$ deleting the vertex $v$ and $v$-shrinking all edges incident with $v$ increases the matching number. The intersecting hypergraphs, as defined above,
are said to be matching critical in the sense that the matching number increases under above
definitions of criticality [M.A. Henning and A. Yeo, Quaest. Math. 37 (2014) 127--138].
Let $n^i(r)$ ($i=2,3,4,5$) denote the  maximum order
of a hypergraph in each  class of matching critical intersecting hypergraphs.
In this paper we study the extremal behavior of matching critical intersecting hypergraphs.  We show that $n^2(r)=n^3(r)$ and $n^4(r)=n^5(r)$ for all $r\ge 2$, which answers an open problem on matching critical intersecting hypergraphs posed by Henning and Yeo. We also give  a strengthening of the result $n^4(r)=n^5(r)$ for intersecting $r$-uniform hypergraphs.

\bigskip

\noindent{\bf AMS (2000) subject classification:}  05C65,  05C70

\vskip 5pt \noindent{\bf Keywords:} Hypergraph; intersecting hypergraph; matching, transversal
\end{abstract}

\section{Introduction}
The relationship between transversals and matchings in hypergraphs have been extensively studied  in \cite{Coc, el, fu, henning3, henning4} and elsewhere. In this paper, we study transversals and matchings in intersecting hypergraphs.

Hypergraphs are a natural generalization of undirected graphs in which ``edges" may consist of more than 2 vertices.
More precisely, a (finite) {\em hypergraph} $H=(V, E)$ consists of a (finite) set $V$ and a collection $E$ of non-empty
subsets of $V$. The elements of $V$ are called {\em vertices} and the elements of $E$ are called {\em hyperedges}, or simply {\em edges} of the
hypergraph. If there is a risk of confusion we will denote the vertex set and the edge set of a hypergraph $H$ explicitly
by $V(H)$ and $E(H)$, respectively.
An $r$-{\em edge} is an edge containing exactly $r$ vertices. The {\em rank} of a hypergraph $H$ is the maximum size of an edge in $H$. Specially, An $r$-{\em uniform} hypergraph $H$ is a hypergraph such that all edges are  $r$-edges.  Obviously,   every (simple) graph is a 2-uniform hypergarph. Throughout this article, all edges have size at least 2.

Two vertices $u$ and $v$ of $H$ are {\em adjacent} in $H=(V,E)$ if there is an edge $e$ in $H$ such that $u,v\in e$. A vertex $v$ and an edge $e$ of $H$ are {\em incident} if $v\in e$. The {\em degree} of a vertex $v \in V$, denoted by $d_{H}(v)$ or  $d(v)$ for shortly if $H$ is clear from
the context, is the number of edges incident to $v$.
The {\em minimum degree} among the vertices of $H$ is denoted by $\delta(H)$.   Two edges in $H$ are said to be {\em overlapping} if they intersect in at least two vertices. The {\em quasidegree} of a vertex $v$ of $H$, denoted by $qd_{H}(v)$ or simply by $qd(v)$, is the maximum number of edges in $H$ whose pairwise intersection is only $v$.
two vertices $u$ and $v$ are {\em connected} if there exists a sequence $u=u_{0}, u_{0}, u_{1}, \ldots, u_{k}=v$ of vertices of $H$ in which $u_{i-1}$ is adjacent to $u_{i}$ for $i=1, 2, \ldots, k$. A {\em connected hypergraph} is a hypergraph in which every pair of vertices are connected.

%For a subset $X\subseteq V(H)$ of vertices in $H$, we define $H-X$ to be the hypergraph obtained from $H$ by deleting the vertices in $X$ and all edges incident with $X$, and deleting resulting isolated vertices, if any.

A set $T\subseteq V$ is called {\em transversal} (also called {\em vertex cover}) of $H=(V,E)$ if it intersects every edge of $H$, i.e., $T\cap e\neq\emptyset$ for all $e\in E$. The minimum
cardinality of the transversals, denoted by $\tau(H)$, is called the {\em transversal number} (also called {\em covering number}) of $H$.
 A subset $M\subseteq E$ is a {\em matching} if every pair of edges from $M$ has an empty intersection. The maximum cardinality
of a matching $M$ is called the {\em matching number}, denoted by $\alpha'(H)$.
 If $M$ is a matching in $H$, then we call a vertex that belongs to an edge of $M$ an $M$-{\em matched vertex}.
Transversals and matchings in hypergraphs have been extensively studied in the literatures (see, for example,
\cite{Alon, chv, Dorf, h, henning1, henning2, k, lai,t}).

%If every vertex has exactly degree one, the matching $\mathcal{M}$ is called perfect. A {\em maximum matching} is a matching $M$ such that $M$ is not a proper subset of another matching.

A hypergraph $H$ is $k$-{\em colorable} if there exists a coloring of the vertices of $H$ using
$k$ colors such that there is no monochromatic edge. A hypergraph $H$ is $k$-{\em chromatic}
if $t$ is the smallest value for which $H$ is $k$-colorable.
%A 2-chromatic hypergraph is also called a {\em bipartite hypergraph}.

%If $e$ is an edge of $H$, then by {\em contracting} the edge $e$ we mean deleting the edge $e$ and replacing it by a smaller, non-empty subset $e'\subset e$. Let $H=(V,E)$ be a hypergraph, and let $e\in E$ and let $v$ be an arbitrary vertex in $e$.
For a subset $E'\subseteq E(H)$ of edges in $H$, we define $H-E$ to be the hypergraph obtained from $H$ by deleting the edges in $X$ and resulting isolated vertices, if any. If $E'=\{e\}$, then we write $H-E'$ simply as $H-e$. If we remove the vertex $v$
from the edge $e$, we say that the resulting edge is obtained by $v$-{\em shrinking} the edge
$e$.

A hypergraph is called {\em intersecting} if any two edges have nonempty intersection. Clearly, $H$ is intersecting if and only if 
$\alpha'(H)=1$. Even for intersecting hypergraphs,
 a long-standing open problem known as Ryser's Conjecture  is open for all $r\ge 6$. Intersecting hypergraphs are well studied in the literature (see, for example, \cite{ah,da, Ekr,fr, gui,hilton,man}).

In \cite{henning4} Henning and Yeo introduced 
 five classes of matching critical intersecting hypergraphs.   In this paper we restrict our attention to intersecting hypergraphs with
the matching criticality. We study the extremal behavior of the matching critical intersecting hypergraphs.
%For example, the famous Erd\"os-Ko-Rado theorem \cite{Ekr} and the Hilton-Milner theorem \cite{hilton}.

It
is clear that if $H$ has rank $r$ then $\tau(H)\le r\alpha'(H)$, and this is attained for example by the
complete $r$-uniform hypergraph $K_{2r-1}^r$ with $2r-1$ vertices, which has $\tau(K^r_{2r-1})=r$ and
$\alpha'(K^r_{2r-2})=1$. In particular, if $H$ is an  intersecting hypergraph with rank $r$, then $\tau(H)\le r$. 
Motivated by this observation, Henning and Yeo \cite{henning4} gave the following definition.

\begin{definition}\label{def1}
For $r\geqslant 2$,  an intersecting hypergraph $H$ of rank $r$ with $\delta(H)\geq 1$  is 1-special if $\tau(H)=r$.
Further, $H$ is maximal 1-special if $H$ is 1-special and adding any missing $r$-edge to $H$ increase the matching number.
\end{definition}

We remark that a 1-special intersecting hypergraph of rank $r$ must be  $r$-uniform, since every edge of $H$ is  a transversal of $H$.

The  two other families of matching criticality in hypergraphs, namely $\alpha'$-edge-criticality and $\alpha'$-vertex-criticality, are  defined in \cite{henning4}.

\begin{definition}\label{def2}
For a hypergraph $H$, $H$ is $\alpha'$-edge-critical if for any $e\in E(H)$ and any $v\in e$, $v$-shrinking $e$ increases the matching number.
$H$ is $\alpha'$-vertex-critical if for every vertex $v$ in $H$ deleting the vertex $v$ and $v$-shrinking all edges incident with $v$ increases the matching number. In particular, when $\alpha'(H)=1$,  $\alpha'$-edge-critical and $\alpha'$-vertex-critical are simply called 1-edge-critical and 1-vertex-critical, respectively.
\end{definition}

By Definition \ref{def2}, we remark that shrinking a 2-edge in a $\alpha'$-edge-critical hypergraph $H$ will yield a 1-edge, although we require
that the original hypergraph $H$ contains no 1-edge. For example, let $H$ be a hypergraph with $V(H)=\{v_1,v_2, v_3, v_4\}$ and $E(H)=\{e_1,e_2,e_3,e_4\}$ where $e_1=\{v_1,v_2, v_3\}$, $e_2=\{v_1,v_4\}$, $e_3=\{v_2,v_4\}$ and $e_4=\{v_3,v_4\}$. It is easy to see that $H$ 1-edge-critical.
From this point of view,  Gallai \cite{ga} showed that the complete graph
$K_{2k+1}$ on $2k + 1$ vertices is the unique $k$-edge-critical connected graph $G$. In particular, $K_3$ is the only 1-edge-critical graph.

Observed by Erd\H{o}s and Lov\'{a}sz \cite{el}, every intersecting hypergraph is 3-colorable. We can get a 3-coloring by coloring the vertices in an arbitrary edge of a intersecting hypergraph with two colors and then coloring all vertices not in this edge with a third color.

For an integer $r\ge 2$, let $\mathcal{H}_{r}$ denote the class of all intersecting hypergraphs of rank $r$ with no  edge consisting of one vertex.
The following five subfamilies of hypergraphs in $\mathcal{H}_r$ are defined in \cite{henning4} (also see \cite{han}).
\begin{eqnarray*}
\mathcal{H}^1(r)&=&\{H\in\mathcal{H}_{r}\mid \mbox{$H$ is 3-chromatic and $r$-uniform}\}.\\
\mathcal{H}^2(r)&=&\{H\in\mathcal{H}_{r}\mid \mbox{$H$ is maximal 1-special}\}.\\
\mathcal{H}^3(r)&=&\{H\in\mathcal{H}_{r}\mid \mbox{$H$ is 1-special}\}.\\
\mathcal{H}^4(r)&=&\{H\in\mathcal{H}_{r}\mid \mbox{$H$ is 1-edge-critical}\}.\\
\mathcal{H}^5(r)&=&\{H\in\mathcal{H}_{r}\mid \mbox{$H$ is 1-vertex-critical}\}.
\end{eqnarray*}
%\begin{definition}\label{def3}
%For $i\in \{1, 2, 3, 4, 5\}$, we define a subclass $\mathcal{H}^{i}_{r}$ of $\mathcal{H}_{r}$ as follows. \\
%$(1)$ A hypergraph $H\in\mathcal{H}^{1}_{r}$ if and only if $H\in\mathcal{H}_{r}$ and $H$ is 3-chromatic and $r$-uniform. \\
%$(2)$ A hypergraph $H\in\mathcal{H}^{2}_{r}$ if and only if $H\in\mathcal{H}_{r}$ and $H$ is maximal 1-special. \\
%$(3)$ A hypergraph $H\in\mathcal{H}^{3}_{r}$ if and only if $H\in\mathcal{H}_{r}$ and $H$ is 1-special. \\
%$(4)$ A hypergraph $H\in\mathcal{H}^{4}_{r}$ if and only if $H\in\mathcal{H}_{r}$ and $H$ is 1-edge-critical. \\
%$(5)$ A hypergraph $H\in\mathcal{H}^{5}_{r}$ if and only if $H\in\mathcal{H}_{r}$ and $H$ is 1-vertex-critical.
%\end{definition}
%Note that if a hypergraph $H$ belongs to $H\in\mathcal{H}^{4}_{r}$ then $v$-shrinking only one edge as opposed to $v$-shrinking all edges incident with $v$ in the definition of $H\in\mathcal{H}^{5}_{r}$. Hence, $H$ is also belong to$H\in\mathcal{H}^{5}_{r}$.
For $r\geqslant2$ and for $i\in \{1, 2, 3, 4, 5\}$, the maximum order of a hypergraph in the class $\mathcal{H}^{i}_{r}$ denoted by $n^{i}(r)$. Thus, $$n^{i}(r)=\max\{|V(H)|\colon\, H\in\mathcal{H}^{i}_{r}\}.$$

 A result gave by Henning and Yeo \cite{henning4} is that the five subclasses of hypergraphs in $\mathcal{H}_{r}$ defined in Definition 4 are nested families, which was also stated in \cite{han} without proof..

\begin{theorem}(\cite{han,henning4}) \label{thm1.1}
For $r\geqslant2$, $\mathcal{H}^{1}_{r}\subset\mathcal{H}^{2}_{r}\subset\mathcal{H}^{3}_{r}\subset\mathcal{H}^{4}_{r}\subset\mathcal{H}^{5}_{r}$.
\end{theorem}

As an immediate consequence of Theorem \ref{thm1.1}, we have the following inequality chain.

\begin{theorem}\label{thm1.2}
(\cite{henning4})
For $r\geqslant2$, $n^{1}(r)\leqslant n^{2}(r)\leqslant n^{3}(r)\leqslant n^{4}(r)\leqslant n^{5}(r)$.
\end{theorem}

 For $r=2,3$,   Henning and Yeo \cite{henning4} proved that the above inequality chain is an equality chain.
\begin{theorem}\label{thm1.3}(\cite{henning4})
For $i\in \{1, 2, 3, 4, 5\}$, $n_i(2)=3$ and $n_i(3)=7$.
\end{theorem}
%\noindent {\bf Open Problem}. For $r\ge 4$, determine whether the inequality chain in Corollary 1 is an equality chain. In particular, is it true that $n_1(r)=n_2(r)=n_3(r)$ for all $r\ge 2$ and is it true that $n_4(r)=n_5(r)$ for all $r\ge 2$?

An open problem posed by Henning and Yeo \cite{henning4} is whether the equalities $n^{1}(r)=n^{2}(r)=n^{3}(r)$  and $n^{4}(r)=n^{5}(r)$  for  $r\geqslant 4$ hold. Further, is it true that the inequality chain in Theorem \ref{thm1.2} is an equality chain for $r\geqslant 4$?

 In this paper we show that $n^{2}(r)=n^{3}(r)$  and $n^{4}(r)=n^{5}(r)$ for all $r\geqslant 2$.
We also prove a strengthening of the  equality $n^{4}(r)=n^{5}(r)$ for intersecting $r$-uniform hypergraphs.

\section{1-Special and maximal 1-special hypergraphs}
In this section we shall prove that $n^{2}(r)=n^{3}(r)$ for all $r\geqslant 2$. To do this,
we first observe the following  relationship between the 1-special hypergraphs and  maximal 1-special hypergraphs.
\begin{proposition} \label{prop2.1}
A hypergraph $H$ is maximal 1-special if and only if $H$ is 1-special and every minimum transversal in $H$ is an edge of $H$.
\end{proposition}

\proof Let $H$ be a maximal 1-special hypergraph. By Definition \ref{def1}, $H$ is 1-special, and so $H$ is $r$-uniform.
 Suppose that  there exists a minimum transversal $T$ of $H$ such that $T$ is not an edge of $H$. Then $|T|=r$.  Hence $T$ can be regarded  as a missing $r$-edge of $H$. Let $H'$ be the hypergraph obtained from $H$ by adding the new $r$-edge $T$ to $H$. Clearly, $\alpha'(H')=\alpha'(H)=1$, contradicting
 the assumption that $H$ is maximal 1-special.

Conversely, suppose that $H$ is a  1-special hypergraph and every minimum transversal in $H$ is an edge of $H$.
We show that $H$ is maximal 1-special. If not, then there exists a missing $r$-edge of $H$ and adding $e$ to $H$ does not increase matching number.
This implies that $e$ has a nonempty intersection with every edge of $H$. Hence $e$ is a minimum transversal in  the resulting hypergraph.
This contradicts the assumption that every minimum transversal in $H$ is an edge of $H$. \endproof

By Proposition \ref{prop2.1}, we now show that $n^{2}(r)=n^{3}(r)$ for all $r\geqslant 2$.
\begin{theorem}\label{thm2.1}
$n^{2}(r)=n^{3}(r)$ for all $r\geqslant 2$.
\end{theorem}

\proof By Theorem \ref{thm1.2}, we have $n^{2}(r)\le n^{3}(r)$ for all $r\geqslant 2$.
To establish the opposite inequality, we show that for
an arbitrary 1-special hypergraph $H$, there exists a maximal 1-special hypergraph $H^*$ such that $V(H)=V(H^*)$.

Let $H$ be an arbitrary 1-special hypergraph, i.e., $H\in \mathcal{H}_r^3$. By Definition \ref{def1}, we have $\tau(H)=r$.
If $H$ is maximal 1-special, there is nothing
to prove. Otherwise,  by Proposition \ref{prop2.1},  there must exist a minimum transversal, say $T$, of $H$ such that $T$ is not an edge of $H$.
We add $T$ to $H$ as a new edge of $H$, and denote by $H'$ the resulting hypergraph. Clearly, $\tau(H')=\tau(H)=r$, so $H'$ is still 1-special.
If $H'$ is  maximal 1-special, then let $H^*=H'$, we are done. Otherwise,  by Proposition \ref{prop2.1} again, there  exists a minimum transversal $T'$ of $H'$ such that $T'$ is not an edge of $H'$.
As above, we obtain a 1-special hypergraph $H''$ by adding $T'$ to $H'$ as a new edge of $H'$.
We repeat the procedure until no minimum transversal $T^*$ in the resulting hypergraph $H^*$ such that $T^*$ is not an edge of $H^*$.
Then $H^*$ is a maximal 1-special hypergraph  such that $V(H)=V(H^*)$, as desired.  Therefore, we have $n^{2}(r)\ge n^{3}(r)$ for all $r\geqslant 2$.
 \endproof

\section{1-Edge-critical and 1-vertex-critical hypergraphs}

In this section we shall show that $n^{4}(r)=n^{5}(r)$ for all $r\geqslant 2$.

By Definition \ref{def2}, we immediately have the following observation.

\begin{observation}\label{ob3.1}
A hypergraph $H$ is 1-edge-critical if and only if
for any $e\in E(H)$ and any $v\in e$, there exists an edge $f\in E(H)$ such that $e\cap f=\{v\}$.
\end{observation}

Henning and Yeo \cite{henning4} provided  a necessary and sufficient condition for a hypergraph to be 1-vertex-critical.

\begin{lemma} (\cite{henning4})\label{lem3.1}
For all $r\geqslant2$, a hypergraph  $H\in H_r^5$
 if and only if $H\in H_r$  and $qd(v)\geqslant2$ for all $v\in V(H)$.
\end{lemma}

\begin{lemma} \label{lem3.2}
For every 1-vertex-critical hypergraph $H\in H^5(r)$, there exists a  1-edge-critical hypergraph $H'\in H^4(r')$ such that
$r'\le r$ and  $V(H)=V(H')$.
\end{lemma}

\proof %Without loss of generality, we may assume that $H$ is minimal 1-vertex-critical.
 There is nothing to prove if  $H$ is 1-edge-critical hypergraph, so we may assume that $H$ is not 1-edge-critical. By Definition \ref{def2}, there exists an edge $e$ of $H$ and a vertex $u\in e$  such that $u$-shrinking $e$ dose not increases the matching number. So  every edge $f$ that contains the vertex $u$ satisfies $|e\cap f|\geqslant2$ by Observation \ref{ob3.1}.
 Since $H\in H^5(r)$, we have $qd(u)\ge 2$ by Lemma \ref{lem3.1}. This implies that
the edge $e$ contribute zero to the quasidegree of $u$, so $u$-shrinking $e$ does not decrease the quasidegree of each vertex in $H$. Now we replace $e$ by the new edge $e_1$ obtained from $e$ by $u$-shrinking $e$. Let $H_1=(H-e)\cup\{e_1\}$.
Then $V(H)=V(H_1)$ and $H_1$ is still 1-vertex-critical, but  $H_1$ possibly has   smaller rank than $H$.
If $H_1$ is not 1-edge-critical, then, by repeating this process of shrinking edges, we obtain a 1-edge-critical hypergraph $H'$ with
$V(H)=V(H')$ and rank $r'$, where $r'\leqslant r$, as desired.
 \endproof

The following lemma is the key to the proof of the main result in this section.
\begin{lemma} \label{lem3.3}
$n^{4}(r+1)>n^{4}(r)$ for all $r\ge2$.
\end{lemma}

\proof Let $H$ be a hypergraph in $\mathcal{H}^{4}_{r}$ with $|V(H)|=n^{4}(r)$. By the definition of 1-edge-critical hypergraphs, we have $d(v)\ge 2$ for each vertex $v\in H$.
Let $e$ be an arbitrary $r$-edge in $H$. We distinguish
the following cases depending on the degree of vertices in $e$.

{\em Case 1:} $d(u)=2$ for all $u\in e$.

 In this case, since $H$ is an intersecting hypergraph and each vertex of $e$ has degree two, $H$ contains exactly $r+1$ edges.
 Let $E(H)=\{e_0, e_1, \ldots, e_r\}$ where $e_0=e$. We construct a 1-edge-critical hypergraph with rank $r+1$ as follows.

 Let $H'$ be the hyergraph obtained from $H$ by
 by adding $r$+1 new vertices $x_{0}, x_{1}, \cdots, x_{r}$ to $H$ and a new edge $\{x_{0}, x_{1}, \cdots, x_{r}\}$ to $H$ and, replacing $e_i$ by $e_i\cup\{x_{i}\}$ for $i=0,1, \ldots, r$.  By our construction, clearly $H'$ is intersecting and it has rank $r+1$. For notational convenience,
let $e'_{r+2}=\{x_{0}, x_{1}, \cdots, x_{r}\}, e_i'=e_i\cup\{x_{i}\}$ for all $i$, $0\le i\le r$. Then $E(H')=\{e'_0, e'_1, \ldots, e'_{r+2}\}.$
We claim that $H'$ is 1-edge-critical hyergraph of rank $r+1$.

By the construction again, we see that $d(x_i)=qd(x_i)=2$ for all $i$, $0\le i\le r$, so  $x_i$-shrinking any edge $e_i'$ containing $x_i$ increases the
matching number.
Therefore, it remains to show that for an arbitrary $e_i'\in E(H')$ ($i\neq r+2$) and $u\in e$ ($u\neq x_i, 0\le i\le r$),
 $u$-shrinking the edge $e_i'$  increases the
matching number. This is indeed so. In fact, since $H$ is 1-edge-critical,  $u$-shrinking the edge $e_i$ will increases the
matching number. Let $\{e_i\setminus \{u\}, e_j\}$ be a matching in the resulting hypergraph when we $u$-shrink the edge $e_i$. Then
clearly $\{e'_i\setminus \{u\}, e'_j\}$ is a matching in the resulting hypergraph when we $u$-shrink the edge $e'_i$.
Hence  $H'\in H^4(r+1)$. Consequently, $n_{r+1}^4\ge |V(H')|=|V(H)|+r+1=n^{4}(r)+r+1>n^{4}(r)$.

{\em Case 2:} There exists a vertex $u\in e$ such that $d(u)\ge 3$.

Since $H\in H^4(r)$, $u$-shrinking the edge $e$ will increase the
matching number, so there exists an edge $f$ such that $f\cap e=\{u\}$.
Let $\{e_0, e_1, \ldots, e_{d(u)-1}\}$ be the set of all edges containing $u$, where $e, f$ are renamed as $e_0, e_{d(u)-1}$, respectively.
We first construct an intersecting hypergraph $H'$ as follows. Let $H'$ be the hypergraph obtained from $H$ by adding two new vertices $\{x, y\}$ and replacing $e$ by $e_0\cup\{x\}$, $e_{d(u)-1}$ by $e_{d(u)-1}\cup\{y\}$, $e_i$ by $(e_{i}\setminus\{u\})\cup\{x, y\}$ for $i=1, 2, \ldots, d(u)-2$.
By the above construction, clearly $H'$ is intersecting and has rank $r+1$.

For notational convenience, we write $e'_0, e'_{d(u)-1}$ and $e'_i$ for $e_0\cup\{x\}$,  $e_{d(u)-1}\cup\{y\}$ and $(e_{i}\setminus\{u\})\cup\{x, y\}$ ($1\le i\le d(u)-2$), respectively.  Let $S=\{e'_0, e'_1, \ldots, e'_{d(u)-1}\}$.
Then $E(H')=(E(H)\setminus\{e_0, e_1, \ldots, e_{d(u)-1}\})\cup S$.

To obtain the desired 1-edge-critical hypergraph. We refine the edges in $S$  by  the following procedure.

{\bf Step 1.} First check every edge $e'_i$ ($0\le i\le d(u)-1$) of $S$ one by one in order. For $e'_i\in S$, if $x$-shrinking $e'_{i}$ dose not increase the matching number, namely $|e'_i\cap e'_j|\ge 2$ for each edge containing $x$ by Observation \ref{ob3.1},  then we replace $e'_{i}$ by $e'_{i}\setminus\{x\}$ and rename $e'_{i}\setminus\{x\}$ as $e'_{i}$.

Let $S'$ (possibly empty) denote the set of the remaining edges $e'_i$ by the above procedure. Then clearly every edge $e'_i$ in $S'$ other than $e'_{d(u)-1}$ still
contains the vertices $x,y$.

{\bf Step 2.} Further check every edge $e'_i$ of $S'$ one by one in order. If $y$-shrinking $e'_{i}$ dose not increase the matching number, namely $|e'_i\cap e'_j|\ge 2$ for each edge containing $y$, then we replace $e'_{i}$ by $e'_{i}\setminus\{y\}$ and rename $e'_{i}\setminus\{y\}$ as $e'_{i}$.

When the procedure terminates, we denote the final resulting hypergraph by $H''$. By the construction, either both or one of $\{x,y\}$ is still in $H''$.
Furthermore, it is easy to see that possibly $qd(x)=0$ (when $e_{i}$ and $e_0$ are overlapping for all $1\leqslant i\leqslant d(u)-1$) or $qd(y)=0$ in $H''$.
Specifically, if $qd(x)=0$, then $qd(y)\ge 2$ in $H''$; if $qd(y)=0$, then $qd(x)\ge 2$ in $H''$.
We show that $H''$ is a 1-edge-critical hypergraph with rank $r+1$.

We first show that $H''$ is intersecting. It suffices to show that every two distinct edges of $\{e'_0, e'_1, \ldots, e'_{d(u)-1}\}$ in $H''$ have
a non-empty intersection since $H$ is intersecting. By the construction of $H''$, $u\in e'_0\cap e'_{d(u)-1}$. We consider every edge
$e'_i$, for $1\le i\le d(u)-2$.
 If $e'_i\cap \{x,y\}=\{y\}$, then  $e'_i$ and each vertex $e'_j$ ($1\le j\le d(u)-2$) containing $y$ share at least the  vertex $y$ in common and, $e'_i$ and each vertex $e'_j$ ($0\le j\le d(u)-2$) containing $x$ (possibly containing $y$) share the common vertex $y$ or intersect within $e'_j$.
If $e'_i\cap \{x,y\}=\{x\}$, then $e'_i$ and each edge $e'_j$ ($0\le j\le d(u)-1$) containing $x$ share at least the  vertex $x$ in common and,
$e'_i$ and each vertex $e'_j$ ($0\le j\le d(u)-2$) containing $y$ (possibly containing $x$) share the common vertex $x$ or intersect within $e'_j$.
If $e'_i\cap \{x,y\}=\{x,y\}$, then $e'_i$ and  $e'_j$ intersect in $x$ or $y$. Thus $H''$ is intersecting.

In order to show that $H''$ is 1-edge-critical, since $H$ is 1-edge-critical, it is enough to show that for any edge $e'_i$ ($0\le i\le d(u)-1$) containing
both or one of $\{x,y\}$ and  any vertex $v\in e'_i\cap \{x,y\}$  in $H''$, $v$-shrinking $e'_i$ can increase the matching number.
Indeed, if $x\in e'_i$, then $x$-shrinking $e'_i$ can increase the matching number, for otherwise we would delete the vertex $x$ from $e'_i$ by Step 1.
If $\{x,y\}\subset e'_i$ or $y\in e'_{d(u)-1}$ , then $y$-shrinking $e'_i$ can increase the matching number, for otherwise we would delete the vertex $y$ from $e'_i$ by Step 2.
If $e'_i\cap \{x,y\}=\{y\}$ ($i\neq d(u)-1$), then there exists an edge $e'_j$ containing $y$ such that $e'_i\cap e'_j=\{y\}$, namely $y$-shrinking $e'_i$ can increase the matching number. Hence, $H''$  is 1-edge-critical.

Note that $H''$ contains at least one vertex in $\{x,y\}$, so $|V(H'')|\ge |V(H)|+1$. Since $H$ is a hypergraph with maximum order in $\mathcal{H}^{4}_{r}$.
This implies that $H''$ has rank $r+1$. Therefore, $n^4(r+1)\ge |V(H'')|>|V(H)|=n^4(r)$.
 \endproof

\begin{theorem}\label{thm3.1}
For any integers $r\ge 2$, $n^{4}(r)=n^{5}(r)$.
\end{theorem}
 \proof By Theorem \ref{thm1.2}, we have $n^4(r)\le n^5(r)$. Let $\mathcal{H}\in H^5(r)$ with $|V(H)|=n^5(r)$. By Lemma \ref{lem3.2},
there exists a  1-edge-critical hypergraph $H'\in H^4(r')$ such that
$r'\le r$ and  $V(H)=V(H')$, so $|V(H')|=n^5(r)$. By Lemma \ref{lem3.3}, we have $n^4(r)\ge n^4(r')\ge |V(H')|=n^5(r)$.
Therefore, $n^{4}(r)=n^{5}(r)$.
 \endproof
\begin{figure}[htbp]
  \centering
  \includegraphics[height=6cm]{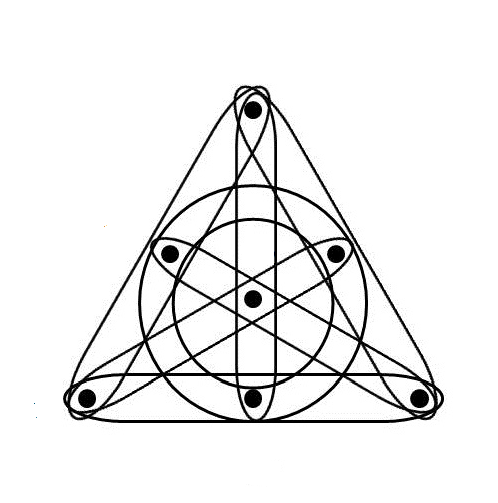}
  \caption{Fano plane: A non-minimal 1-vertex-critical hypergraph}\label{Fig.1}
\end{figure}
Now we give a strengthening of Theorem \ref{thm3.1}. We shall show that the maximum order $n^5(r)$ of 1-vertex-critical hypergraphs
is the same as that of 1-vertex-critical $r$-uniform hypergraphs.

We call a 1-vertex-critical hypergraph $H$
is {\em minimal} if $H$ contains no edge $e$ such that
$H-e$ is 1-vertex-critical. For $r\geqslant2$, let  $\mathcal{H}_{uni}^5(r)$ and $\mathcal{H}_{min}^5(r)$  be the sets of1-vertex-critical $r$-uniform
 and minimal 1-vertex-critical
  hypergraphs in $\mathcal{H}^5(r)$, respectively.
 Obviously,
$\mathcal{H}_{uni}^5(r)\subseteq \mathcal{H}^5(r)$, $\mathcal{H}_{min}^5(r)\subseteq \mathcal{H}^5(r)$, but $\mathcal{H}_{uni}^5(r)\nsubseteq \mathcal{H}_{min}^5(r)$. Fano plane  in Fig. 1 clearly is a non-minimal 1-vertex-critical hypergraph.
Let $n^5_r={\rm max}\{|V(H)|\colon\, H\in \mathcal{H}_{uni}^5(r)\}$ and $n^4_r={\rm max}\{|V(H)|\colon\, H\in \mathcal{H}^4(r)~ \mbox{and $H$ is $r$-uniform}\}$.

\begin{lemma}\label{lem3.4}
For every $H\in \mathcal{H}^5(r)$, there exists an $H'\in H_{min}^5(r)$  such that $V(H)=V(H')$, and every edge $e\in E(H')$ contains at least a vertex $v$ such that $qd_{H'-e}(v)=0$.
%$d(u)=qd(u)=2$ or $qd(u)=2$ and every pair of all edges  which contain $u$ but $e$ are overlapping.
\end{lemma}
\proof
If $H$ is minimal 1-vertex-critical, we are done.  Otherwise, there exists an edge $e_1\in E(H)$ such that $qd_{H-e_1}(v)\ge 2$  for all $v\in V(H-e_1)$ in $H-e_1$. We delete the edge $e_1$ from $H$.
By Lemma \ref{lem3.1},  the  hypergraph $H-e_1$ is still 1-vertex-critical and $V(H)=V(H-e_1)$.
If $H-e_1$ still contains an edge $e_2$ such that $qd_{H-\{e_1,e_2\}}(v)\ge 2$ for all $v\in V(H-\{e_1,e_2\})$, then
we delete the edge $e_2$ from $H-e_1$. Then, by Lemma \ref{lem3.1}, $H-\{e_1,e_2\}$ is still 1-vertex-critical.
By repeating this above process of deleting edges in $H$, we obtain
a hypergraph $H'$ still 1-vertex-critical but $H'$ contains no edge $e$ such that
$H'-e$ is 1-vertex-critical. Therefore, $H'\in H_{min}^5(r)$ and $V(H)=V(H')$.
\endproof

By Lemma \ref{lem3.4}, it is easy to see that $n^5(r)={\rm max}\{|V(H)|\colon\, H\in \mathcal{H}_{min}^5(r)\}$. Furthermore, we have
the following lemma.

\begin{lemma} \label{lem3.5}
For $r\geqslant2$, every hypergraph $H$ in $H_{min}^5(r)$ with $|V(H)|=n^5(r)$  is $r$-uniform.
\end{lemma}

\proof Suppose, to the contrary, that $H$ contains a $t$-edge $e$ such that $t<r$.
 Since $H$ is minimal 1-vertex-critical,  there exists a vertex $u$ in $e$ such that $qd_{H-e}(u)<2$ by Lemma \ref{lem3.4},
 hence there exists an edge $f$  of $H$ whose intersecting with $e$ is only $u$, i.e., $e\cap f=\{u\}$.  Clearly, $f\setminus \{u\}$ is a transversal in $H-e$.
 Note that $H$ has rank $r$,
without loss of generality, we may assume that $f$ is an $r$-edge of $H$.

Let $H'$ be the hypergraph obtained from $H$ by adding a new vertex $v$ and a new edge $(f\setminus \{u\})\cup\{v\}$ to $H$ and replacing $e$ by $e\cup\{v\}$.
 By the construction of $H'$, $V(H')=V(H)\cup\{v\}$ and $E(H')=(E(H)\setminus \{e\})\cup \{e\cup \{v\}, (f\setminus \{u\})\cup\{v\}\}$.
Hence $|V(H')|=n^5(r)+1$ and $H'$ has rank $r$. Note that $f\setminus \{u\}$ is a transversal of $H-e$, so $(f\setminus \{u\})\cup\{v\}$ meets with all edges of $H'$.
Hence $H'$ is a intersecting hypergraph with rank $r$. On the other hand, we note that $d(v)=qd(v)=2$ and $qd(u)\geqslant2$ for all $u\in H'$.
Thus $H'$ is still 1-vertex-critical, i.e., $H'\in \mathcal{H}^5(r)$. But then $|V(H')|\le n^5(r)$, which is a contradiction. \endproof

By  Lemma \ref{lem3.5} and Theorem \ref{thm3.1}, we obtain a strengthening of Theorem \ref{thm3.1}.
\begin{theorem}\label{thm3.2} For any integers $r\ge 2$, $n^4_r=n^4(r)=n^5(r)=n^5_r$.
\end{theorem}
\proof
By Lemma \ref{lem3.5}, we have $n^5_r=n^5(r)$. Let $H\in\mathcal{H}_{min}^5(r)$ is a hypergraph with order $n^5_r$. We claim that $H$ is 1-edge-critical.
If not, then there exists an edge $e$ and a vertex $v\in e$ such that replacing $e$ by $e'$ which obtained from $v$-shrinking $e$ dose not increase the matching number. It implies that $e$ contribute zero to the quasidegree of $v$. Hence, there exists a hypergraph $H'\in H^5(r)$ with $|V(H')|=n^5_r$ such that $H'$ there exists an $(r-1)$-edge $e'$.
But this  contradicts  Lemma \ref{lem3.5}. Thus $n^5(r)=n^5_r\le n^4_r\le n^4(r)$. By Theorem \ref{thm3.1},
$n^4_r=n^4(r)=n^5(r)=n^5_r$.
\endproof
%²Î¿¼ÎÄÏ×£º°´×Öĸ˳Ðò

%\end{CJK*}


\begin{thebibliography}{99}
\bibitem{ah} Ron Aharoni, C.J. Argue, Covers in partitioned intersecting hypergraphs, European J. Combin. 51 (2016) 222--226.

\bibitem{Alon} N. Alon,  Transversal numbers of uniform hypergraphs, Graphs Combin. 6 (1990) 1--4.

\bibitem{chv} V. Chv\'{a}tal, C. McDiarmid, Small transversals in hypergraphs, Combinatorica 12 (1992) 19--26.

\bibitem{Coc} E.J. Cockayne, S. T. Hedetniemi, P. J. Slater, Matchings and transversals in hypergraphs, domination and independence-in trees, J. Combin. Theory Ser. B 27 (1979) 78--80.
    
\bibitem{da} S. Das, B. Sudakov,  Most probably intersecting hypergraphs, Electron. J. Combin.  22 (1) (2015)\#P 1.80.

\bibitem{Dorf} M. Dorfling, M. A. Henning, Linear hypergraphs with large transversal number and maximum degree two, European J. Combin. 36 (2014) 231--236.

\bibitem{Ekr} P. Erd\"os, C. Ko and R. Rado, Intersection theorems for systems of finite sets,
Quart. J. Math. Oxford 12(2) (1961) 313--320.

\bibitem{el} P. Erod\"{o}s, L. Lov\'{a}sz, Problems and results on 3-chromatic hypergraphs and some related questions. Infinite and finit sets (Colloq. , Keszthely, 1973; dedicated to P.  Erd\"{o}s on his 60th birthday), Vol. II, pp. 609--627. Colloq Math. Soc. Janos Bolyai 10 North-Holland, Amsterdam, 1975.

\bibitem{ga}T. Gallai, Neuer Beweis eines Tutteschen Satzes, Magyar Tud. Akad. Mat. Kut. Int.
K\"ozl. 8 (1963) 135--139.

\bibitem{fr} P. Frankl, Z. F\"uredi,  Finite projective spaces and intersecting hypergraphs, Combinatorica  6  (1986) 335--354.

\bibitem{fu}Z. F\"uredi, Matchings and covers in hypergraphs, Graphs Combin. 4 (1988) 115--206.

\bibitem{gui}B. Guiduli, Z. Kirfily, On intersecting hypergraphs, Discrete Math. 182 (1998) 139--151.

\bibitem{han}D. Hanson and B. Toft, On the maximum number of vertices in $n$-uniform cliques,
Ars. Combin. 16A (1983) 205--216.

\bibitem{h} P.E. Haxell, A condition for matchability in hypergraphs, Graphs.  Combin. 11 (1995) 245--248.

\bibitem{henning1} M.A. Henning, C. L\"{o}wenstein, Hypergraphs with large transversal number and with edge sizes at least four, Cent. Eur. J. Math.  10(3) (2012) 1133--1140.

\bibitem{henning2} M.A. Henning, A. Yeo, Hypergraphs with large transversal number and with edge sizes at least three, J. Graph Theory 59 (2008) 326--348.

\bibitem{henning3} M.A. Henning, A. Yeo, Transversals and matchings in $3$-uniform hypergraphs, European J. Combin. 34 (2013) 217--228.

\bibitem{henning4} M.A. Henning, A. Yeo, Matching critical intersecting hypergraphs, Quest. Math. 37 (2014), 127--138.

\bibitem{hilton}A.J.W. Hilton and E.C. Milner, Some intersection theorems for systems of finite
sets, Quart. J. Math. Oxford 18(2) (1967) 369--384.

\bibitem{k} D. K$\ddot{u}$hn D. Osthus, Matchings in hypergraphs of large minimum degree, J. Graph Theory 51 (2006) 269--280.

\bibitem{lai} F.C. Lai, G.J. Chang, An upper bound for the transversal numbers of $4$-uniform hypergraphs, J. Combin. Theory Ser. B 50 (1990) 129--133.

\bibitem{man}T. Mansour, C. Song, R. Yuster, A comment on Ryser's conjecture for intersecting
hypergraphs, Graphs  Combin. 25 (2009) 101--109.



\bibitem{t} Zs. Tuza, Critical hypergraphs and intersecting set-pair systems, J. Combin. Theory Ser. B 39 (1985) 134--145.
\end{thebibliography}
\end{document}